\documentclass[12pt,pdftex]{article}
\usepackage[dvips]{graphicx, color}            %
\usepackage{pslatex}	    %

\usepackage{amssymb, amscd, amstext}
\usepackage{latexsym}
\usepackage{amsxtra, amsfonts}

\def\part#1{\frac{\partial\phantom{q}}{\partial#1}}

\makeatletter
 \newlength{\typesize}
\makeatother

\newlength{\vvoff}
\newlength{\hhoff}

\newcommand{\locateoffcenter}[1]{%
\addtolength{\vvoff}{-0.25\typesize}%
\raisebox{\vvoff}{\hspace{\hhoff}\makebox(0,0){\smash{#1}}}
}
\newcommand{\object}[1]{%
\setlength{\vvoff}{0pt}%
\setlength{\hhoff}{0pt}%
\locateoffcenter{#1}
}
\newcommand{\elabel}[1]{%
\setlength{\vvoff}{0.75\typesize}%
\setlength{\hhoff}{0pt}%
\locateoffcenter{#1}
}

\newcommand{\slabel}[1]{%
\setlength{\vvoff}{0pt}%
\setlength{\hhoff}{0.75\typesize}%
\locateoffcenter{#1}
}

\newcommand{\swlabel}[1]{%
\setlength{\vvoff}{-0.5\typesize}%
\setlength{\hhoff}{0.75\typesize}%
\locateoffcenter{#1}
}

\newcommand{\selabel}[1]{%
\setlength{\vvoff}{0.5\typesize}%
\setlength{\hhoff}{0.75\typesize}%
\locateoffcenter{#1}
}

\newcommand\beq{\begin{equation}}
\newcommand\eeq{\end{equation}}

\newcommand{\HH}{\text{\rm H}}

\newcommand{\IF}{\mathbb{F}}

\newcommand{\IP}{\mathbb{P}}                                     
                           
\newcommand{\IR}{\mathbb{R}}                           
\newcommand{\IC}{\mathbb{C}}

\newcommand{\M}{\mathcal{M}}

\newcommand{\cC}{\mathcal{C}}

\newcommand{\cF}{\mathcal{F}}

\newcommand{\g}{       \mathfrak{g}     }

\newcommand{\gl}{       \mathfrak{gl}     } %
\newcommand{\lsl}{       \mathfrak{sl}     } %

\newcommand{\PVI}{{$\text{\rm P}_{\text{\rm VI}}$}}   %
\newcommand{\pf}{\begin{bpf}}

\newcommand{\pfms}{\begin{bpfms}}
\newcommand{\epf}{\end{bpf}\hfill$\square$\\}           %
\newcommand{\epfms}{\end{bpfms}\hfill$\square$\\}               %

\newcommand{\idea}{\begin{bidea}}

\newcommand{\eidea}{\end{bidea}\hfill$\square$\\}           %

\newcommand{\sk}{\begin{bsk}}    %

\newcommand{\esk}{\end{bsk}\hfill$\square$\\}           %
\newcommand{\sketch}{\begin{bsketch}}%

\newcommand{\esketch}{\end{bsketch}\hfill$\square$\\}

\newcommand{\al}{\alpha}

\newcommand{\be}{\beta}
\newcommand{\ga}{\gamma}
\newcommand{\de}{\delta}

\newcommand{\la}{\lambda}

\newcommand{\tr}{\mathop{\rm Tr}}

\newcommand{\Hom}{\mathop{\rm Hom}}

\newcommand{\SL}{{\mathop{\rm SL}}}
\newcommand{\PSL}{{\mathop{\rm PSL}}}
\newcommand{\GL}{{\mathop{\rm GL}}}
\newcommand{\PGL}{{\mathop{\rm PGL}}}

\newcommand{\SO}{{\mathop{\rm SO}}}

\newcommand{\diag}{{\mathop{\rm diag}}}

\newtheorem {thm}{Theorem}

\begin{document}
\date{}

\title{Towards a nonlinear Schwarz's list} 
\author{Philip Boalch}

\maketitle

\section{Introduction}
 
The main theme of this talk is ``icosahedral" solutions of (ordinary)
differential equations, a topic that seems suitable for a 60th birthday conference. We will however try to go beyond the icosahedron, to see what comes next, and consider various symmetry groups each of which could be thought of as the next in a sequence, following the icosahedral group.

To fix ideas let us give a classical example. 
Recall the icosahedral rotation group of order $60$: 
$$A_5 \cong \PSL_2(\IF_5)\cong \Delta_{235} \cong
\langle\  a,b,c\ \bigl\vert\ a^2=b^3=c^5=abc=1\ \rangle.$$

This is described via three generators $a,b,c$ whose product is the identity, and so it is natural to look for ordinary differential equations on the 3 punctured sphere 
$\IP^1(\IC)\setminus\{0,1,\infty\}$ 
with monodromy group $A_5$. 
Now $A_5$ is a three-dimensional rotation group so naturally lives in $\SO_3(\IR)$ which is a subgroup of $\SO_3(\IC)$ which is isomorphic to $\PSL_2(\IC)$.
Thus we are led to search for connections 
\beq\label{eq: hg con}
 \nabla= d-\left(\frac{A_1}{z}+\frac{A_2}{z-1}\right)dz,\qquad A_i\in\lsl_2(\IC) 
\eeq
on rank two holomorphic vector bundles over the three-punctured sphere with {\em projective} monodromy group equal to $A_5$.

Such connections are essentially the same as Gauss hypergeometric equations, and H. Schwarz \cite{Schwarz} classified all such equations having finite monodromy groups in 1873.  The list he produced has $15$ rows, one for the family of dihedral groups, two rows for each of the tetrahedral and octahedral groups and $10$ rows for the icosahedral group. See Table 1.

\begin{table}[h]
	\resizebox{\textwidth}{!}{\includegraphics[width=\textwidth]{schwarz4small}}
	\caption{Schwarz's list \cite{Schwarz}} \end{table}

A key point here is that the Gauss hypergeometric equation is {\em rigid} so the full monodromy representation (of the fundamental group of the $3$-punctured sphere into $\PSL_2(\IC)$) is determined by the conjugacy classes of the monodromy around each of the  punctures. Thus in Schwarz's list it is sufficient to list these local monodromy conjugacy classes in order to specify the
possible monodromy representations (and from this it is easy to find a hypergeometric equation with given monodromy).
To ease recognition, to the left of the table we have listed the triples of conjugacy classes which occur, labelling the four nontrivial conjugacy classes of $A_5$ by $a,b,c,d$, representing rotations by 
$\frac12,\frac13,\frac15,\frac25$ of a turn, respectively. (In the octahedral case one may also have rotations by a quarter of a turn, which we label by $g$.)

\subsection{Naive generalisations}

Our basic aim is to discuss three naive generalisations of Schwarz's list, as follows.
The first two arise simply by looking for {\em nonrigid} connections that are natural generalisations of the hypergeometric connections considered above, obtained by adding an extra singularity---the two cases are generalisations of two ways one may view the hypergeometric equation as a connection.
First of all we can simply add another pole at some point $t$:
$$
{\bf (A)}  \qquad\qquad 
\nabla= d-\left(\frac{A_1}{z}+\frac{A_2}{z-t}+\frac{A_3}{z-1}\right)dz,\qquad A_i\in\lsl_2(\IC) $$
and keep the coefficients in $\lsl_2(\IC)$.  

Secondly we recall that 
the connection one obtains immediately upon choosing a cyclic vector for the hypergeometric equation is as in 
(\ref{eq: hg con}) but with $A_1,A_2$ both rank one matrices (in $\gl_2(\IC)$). Then the monodromy group will be a complex reflection group (generated by two two-dimensional complex reflections\footnote{i.e. arbitrary automorphisms of the form ``one plus rank one", not necessarily of order two or orthogonal.}) and the natural generalisation is then to consider connections of the form:
$$
{\bf (B)}  \qquad\qquad 
\nabla= d-\left(\frac{B_1}{z}+\frac{B_2}{z-t}+\frac{B_3}{z-1}\right)dz,\qquad B_i\in\gl_3(\IC)$$
with each $B_i$ having rank one, so the monodromy group will be generated by three three-dimensional complex reflections.
This is a very natural condition as we will see.

{\bf Questions A, B:} \ Find the analogue of Schwarz's list for connections {\bf (A)} or {\bf (B)}.

These questions can now be answered and lead to two ``nonrigid Schwarz's lists", i.e. to classifications of possible monodromy representations with finite image (up to equivalence) and the construction of connections realizing such representations.  We should emphasise that the main focus has been the construction of such connections with given monodromy representation for any value of $t$ (which is a tricky business in this nonrigid case), rather than just the classification (which is reasonably straightforward).

{\bf Example} (of type {\bf (B)}).
The full symmetry group of the icosahedron is the icosahedral reflection group of order 120:
$$
H=H_3\cong \langle \ r_1, r_2, r_3\ \bigl\vert\ r_i^2=1, 
(r_1r_2)^2=(r_2r_3)^3=(r_3r_1)^5=1 \ \rangle
$$$$
\subset O_3(\IR)\subset {\GL}_3(\IC). 
$$
This is generated by three reflections (whose product is {\em not} the identity) and so it is natural to look for connections on rank three bundles over a four-punctured sphere with monodromy $H$ (generated by 3 reflections about 3 of the punctures---i.e. connections of the form {\bf (B)} with each of the three residues $B_i$ having trace $\frac{1}{2}$ so the corresponding reflections are of order two). There turn out to be three inequivalent triples of generating reflections of $H$, two of which are related by an outer automorphism. The problem is to write down connections 
of the desired form for any value $t$ of the final pole position. 
One triple of generating reflections is intimately related to K. Saito's flat structure for $H$ (or icosahedral Frobenius manifold)
\cite{KSaito1} 
and appears in Dubrovin's article \cite{Dub95long} Appendix E.
The other two triples were dealt with around 1997 by Dubrovin and Mazzocco \cite{DubMaz00}; one is similar to the first case (since related to it by an outer automorphism) but the final triple turned out to be much trickier, and writing out the family of connections in this case involved a specific elliptic curve which took about ten pages of 40 digits integers to write down 
(see the preprint version of \cite{DubMaz00} on the mathematics arxiv).
We will eventually see below that  this elliptic solution is in fact equivalent to a solution with a simple parametrisation, agreeing with Hitchin's philosophy that ``nice problems should have nice solutions".

{
\renewcommand{\baselinestretch}{0.93}              %
\small
{\bf Remark\ } Before moving on to the third generalisation let us add some other historical comments.
The ``non-naive" generalisations  of the Gauss hypergeometric equation are the equations satisfied by the $\!\!\!\ _nF_{n-1}$ hypergeometric functions (the Gauss case being that of $n=2$). 
The corresponding Schwarz's 
list appears in the 1989 article \cite{BH89} of Beukers and Heckman. In terms of connections this amounts to considering 
connections (\ref{eq: hg con}) on rank $n$ vector bundles, still with three singularities on $\IP^1$, but with $A_1$ of rank $n-1$ and $A_2$ of rank $1$; {\em these connections are still rigid}.
\newline
{\phantom{123}Some work in the nonrigid case has been done (besides that we will recall below) by considering generalisations of the hypergeometric equation as an equation (rather than as a connection); for example the algebraic solutions of the Lam\'e equation were studied in \cite{B-vdW} by Beukers and van der Waall (Lam\'e equations are basically the second order Fuchsian equations with four singular points on $\IP^1$ such that three of the local monodromies are of order two). In general connections of type {\bf (A)} with such monodromy representations will not come from a Lam\'e equation (since upon choosing a cyclic vector the corresponding equations will in general have additional apparent singularities; this can also be seen by counting dimensions). Indeed it turns out (\cite{B-vdW}) that Lam\'e equations only have finite monodromy for special configurations of the four poles.

\renewcommand{\baselinestretch}{1}
\small
}
\renewcommand{\baselinestretch}{1}
\normalsize
\vspace{-0.3cm}
\subsection{Nonlinear analogue---the Painlev\'e VI equation}
One reason hypergeometric equations are interesting is that they provide the simplest explicit examples of {\em Gauss--Manin connections}. Indeed this is one reason  Gauss was interested in them: he observed that the periods of a family of elliptic curves satisfy a (Gauss) hypergeometric equation.
(The modern interpretation of this is as the explicit form of the natural flat connection on the vector bundle of first cohomologies over the base of the family of elliptic curves, written with respect to the basis given by the holomorphic one forms---and their derivatives---on the fibres.) Nowadays there is still much interest in such linear differential equations ``coming from geometry".

Thus the nonlinear analogue of the Gauss hypergeometric equation should be the explicit form of the simplest {\em nonabelian} Gauss--Manin connection (i.e the explicit form of the natural connection on the bundle of first nonabelian cohomologies of some family of varieties). The simplest interesting case corresponds to taking  the universal family of four punctured spheres  and taking cohomology with coefficients in $\SL_2(\IC)$ (one needs a non-trivial family of varieties with nonabelian fundamental groups).
This leads to the Painlev\'e VI equation (\PVI), which is a second order nonlinear differential equation whose solutions, like those of the hypergeometric equation, branch only at $0,1,\infty\in \IP^1$.
In particular we may study the (nonlinear) monodromy of solutions of \PVI, by examining how solutions vary upon analytic continuation along paths in the three-punctured sphere.

Thus, since Schwarz lists fundamental solutions of hypergeometric equations having finite monodromy, our main question is to construct the analogue of Schwarz's list for \PVI:

{\bf Question C:}
What are the solutions of Painlev\'e VI having finite monodromy?

This question is still open; there is as yet no full classification---the main effort (at least of the present author) has been towards finding and constructing interesting solutions. So far all known finite branching solutions are actually algebraic\footnote{Apparently Iwasaki \cite{iwasaki-fbr} has recently shown that all finite branching solutions are algebraic.}. Currently we are at the reasonably happy state of affairs that all such solutions known to exist have actually been constructed. In what follows I will explain various aspects of the problem, and in particular show how the nonrigid lists {\bf (A)} or {\bf (B)} map to the list of {\bf (C)}.
Some key points, demonstrating the richness and variety of solutions, are:

$\bullet$ There are algebraic solutions of \PVI\  not related to finite subgroups of the coefficient group $\SL_2(\IC)$,

$\bullet$ There are `generic' solutions of \PVI\  with finite monodromy; i.e. not lying on any of the reflection hyperplanes of the affine $F_4$ Weyl group of symmetries of  \PVI,

$\bullet$ There are entries on the list of {\bf (C)} which do not come from either {\bf (A)} or {\bf (B)}.

In particular we will see \PVI\  solutions related to the groups
$A_6$, $\PSL_2(\IF_7)$ and $\Delta_{237}$.

\section{What is Painlev\'e VI?}

There are various viewpoints, and simply giving the explicit equation is perhaps the least helpful introduction to it. In brief, Painlev\'e VI is:

$\bullet$ the explicit form of the simplest nonabelian Gauss--Manin connection,
 
$\bullet$ the equation controlling the ``isomonodromic deformations" of certain logarithmic connections/Fuchsian systems on $\IP^1$,
 
$\bullet$ the most general second order ODE with the so-called `Painlev\'e property',
 
$\bullet$ a certain dimensional reduction of the anti-self-dual Yang--Mills equations (see e.g. \cite{MasWoo}),
 
$\bullet$ an equation related to certain elliptic integrals with moving endpoints (cf. R. Fuchs \cite{RFuchs} and Manin \cite{Manin-P6}),

$\bullet$ The second order ODE for a complex function $y(t)$

$$\frac{d^2y}{dt^2}=
\frac{1}{2}\left(\frac{1}{y}+\frac{1}{y-1}+\frac{1}{y-t}\right)
\left(\frac{dy}{dt}\right)^2
-\left(\frac{1}{t}+\frac{1}{t-1}+\frac{1}{y-t}\right)\frac{dy}{dt}  $$$$
\quad\ +\frac{y(y-1)(y-t)}{t^2(t-1)^2}\left(
\al+\be\frac{t}{y^2} + \gamma\frac{(t-1)}{(y-1)^2}+
\delta\frac{ t(t-1)}{(y-t)^2}\right) $$
where $\al,\be,\ga,\de\in\IC$ are constants.

The Painlev\'e property means that any local solution $y(t)$ defined in a disc in the three-punctured sphere 
$\IP^1\setminus\{0,1,\infty\}$ extends to a meromorphic function on the universal cover of $\IP^1\setminus\{0,1,\infty\}$. It is this property that enables us to speak of the monodromy of \PVI\  solutions. Concerning solutions there is a basic trichotomy (see Watanabe \cite{watanabePVI}):
$$
\text{A solution of \PVI\  is either\ \ }
\begin{cases}
\text{a `new' transcendental function, or}\\
\text{a solution of a 1st order Riccati equation, or}\\
\text{an algebraic function.}
\end{cases}
$$
In particular if one is interested in constructing new explicit solutions of \PVI\ then, since the Riccati solutions are all well understood,
the algebraic solutions are the first place to look.

The standard approach to \PVI\ is as isomonodromic deformations of rank two logarithmic connections with four poles on $\IP^1,$ as the poles move (generic such connections are of the form {\bf (A)}, and then $t$ parametrises the possible pole configurations). In particular one can see the four constants in \PVI\ directly in terms of the eigenvalues of the residues of the connection:
if we set $\theta_i$ to be the difference of the eigenvalues (in some order) of the residue $A_i$ ($i=1,2,3,4$, where $A_4=-\sum_1^3A_i$ is the residue at infinity) then the relation to the constants is
$$\al=(\theta_4-1)^2/2, \qquad \be=-\theta_1^2/2, \qquad 
\ga=\theta^2_3/2, \qquad \delta=(1-\theta_2^2)/2.
$$
Before going into more detail let us also mention one further property of \PVI: it admits a group of symmetries isomorphic to the affine Weyl group of type $F_4$ (see \cite{OkaPVI} or the exposition in \cite{srops}). 
Indeed
treating $\theta=(\theta_1,\ldots,\theta_4)\in\IC^4$ as the set of parameters for \PVI\ is useful since the affine $F_4$ Weyl group of symmetries acts in the standard way on this $\IC^4$.
(We will see below that these four parameters may also be interpreted as coordinates on the moduli space of cubic surfaces.)

\subsection{Conceptual approach to $\text{\bf P}_{\text{\bf VI}}$}

Consider the universal family of smooth four-punctured rational curves with labelled punctures.
Write $B := \M_{0,4}\cong\IP^1\setminus\{0,1,\infty\}$ for the base, $\cF$ for the standard fibre and $\cC$ for the total space:

\[
\begin{CD}
\cC @<<< \cF \cong \IP^1\setminus\text{4 points}\\
@VVV\\
B 
\end{CD}
\]

Now replace each fibre $\cF$ by $\HH^1(\cF,G)$ where $G=\SL_2(\IC)$.
Here we will use two viewpoints/realisations of this nonabelian cohomology set $\HH^1$:

1) (Betti): Moduli of fundamental group representations:

$$\HH^1(\cF,G)\cong \Hom(\pi_1(\cF),G)/G$$

2) (De Rham): Moduli of connections on holomorphic vector bundles over $\cF$

These two viewpoints are related by the Riemann--Hilbert correspondence (the nonabelian De Rham functor), taking connections to their monodromy representations.
The point is that algebraically these realisations of $\HH^1$ are very different and the Riemann--Hilbert map is transcendental (things written in algebraic coordinates on one side will look a lot more complicated from the other side).

Thus we get two nonlinear fibrations over the base $B$, with fibres the De Rham or Betti realisations of $\HH^1(\cF,G)$ respectively:

\[
\begin{CD}
\M_{\text{De Rham}} @>\text{Riemann--Hilbert}>> \M_{\text{Betti}}\\
@VVV @VVV\\
B @= B 
\end{CD}
\]

As in the case with abelian coefficients one still gets a natural connection on these cohomology bundles. The surprising fact is that it is algebraic on {\em both} sides (approximating the  De Rham side in terms of logarithmic connections to give it an algebraic structure \cite{Nit-log}).
Thus when written explicitly we will get {\em nonlinear  algebraic} differential equations ``coming from geometry". (See Simpson \cite{Sim94ab} Section 8 for more on these connections in the case of families of projective varieties.)

The two standard descriptions of the abelian Gauss--Manin connection generalise to descriptions of this nonlinear connection. In the Betti picture we may identify two nearby fibres of $\M_{\text{Betti}}$
simply by keeping the monodromy representations (points of the fibres) constant: moving around in $B$ amounts to deforming the configuration of four points in $\IP^1$ and it is easy to see how to identify the fundamental groups of the four punctured spheres as the punctures are deformed---use the same generating loops. This `isomonodromic' description, preserving the monodromy representation, is the nonabelian analogue of keeping the periods of one-forms constant.

On the De Rham side the nonlinear connection can be described in terms of extending a connection on a vector bundle over a fibre $\cF$, to a flat connection on a vector bundle over a family of fibres and then restricting to another fibre, much as the abelian case is described in terms of closed one-forms (linear connections replacing one-forms and flatness replacing the notion of closedness).

Each of these descriptions has a use: the De Rham viewpoint lends itself to giving an explicit form of the nonlinear connection (essentially amounting to the condition for the flatness of the connection over the family of fibres). The Betti viewpoint is more global and allows us to study the monodromy of the nonlinear connection, as an explicit action on fibres of $\M_{\text{Betti}}$.

\subsection{Explicit nonlinear equations}

The De Rham bundle $\M_{\text{De Rham}}$ is well approximated by the space of logarithmic connections with four poles on the trivial rank two holomorphic bundle (with trivial determinant) over $\IP^1$. Call the space of such connections $\M^*$ and observe it parametrises connections of the form {\bf (A)}, and that these are determined by the value of $t\in B$ and the residues:

$$\M^* \cong B\times \left\{ (A_1,\ldots,A_4)\ \bigl\vert\ 
A_i\in \g, \sum A_i=0 \right\}/G.$$
Here $G=\SL_2(\IC)$ does not act on $B$ and acts by diagonal conjugation on the residues $A_i$.
In general this quotient will not be well behaved, but it has a natural Poisson structure and the generic symplectic leaves will be smooth complex symplectic surfaces. Clearly $\M^*$ is trivial as a bundle over $B$ (projecting onto the configuration of poles), but the nonabelian Gauss--Manin connection is different to this trivial connection and was computed about 100 years ago by Schlesinger (essentially in the way stated above it seems). 
The nonlinear connection is given by {\em Schlesinger's equations}, which in the case at hand are:
$$\frac{dA_1}{dt}=\frac{[A_2,A_1]}{t},\qquad 
\frac{dA_3}{dt}=\frac{[A_2,A_3]}{t-1}$$
together with a third equation for $dA_2/dt$ easily deduced from the fact that $A_4$ remains constant. If the residues of the connection satisfy these equations then the corresponding monodromy representation remains constant as $t$ varies. (They are easily derived from the vanishing of the curvature of the `full' connection $d-\left(A_1\frac{dz}{z}+A_2\frac{dz-dt}{z-t}+A_3\frac{dz}{z-1}\right)$.)

To get from here to \PVI\  one chooses specific functions $x,y$ on $\M^*$ which restrict to coordinates on each generic symplectic leaf and writes down the connection in these (carefully chosen) coordinates. (See \cite{k2p} pp.199-200 for a discussion of the formulae, which are from \cite{JM81}.) This leads to two coupled nonlinear first order equations, and eliminating $x$ leads to the second order 
Painlev\'e VI equation for $y(t)$. It was first written down in full generality by R. Fuchs \cite{RFuchs} (whose father L. Fuchs was also the father `Fuchsian equations').

\subsection{Monodromy of Painlev\'e VI}

Since the Betti and De Rham realisations are isomorphic, we see the monodromy of solutions to \PVI\ thus corresponds to the monodromy of the connection on $\M_{\text{Betti}}$. This amounts to an {\em action} of the fundamental group of the base $B$ on a fibre, and this action can be described explicitly.

Let $\M_t=\Hom(\pi_1(\IP^1\setminus\{0,t,1,\infty\}),G)/G$ be the fibre of $\M_{\text{Betti}}$ at some fixed point $t\in B$.
The key point is that $\pi_1(B)\cong \cF_2$ (the free nonabelian group on two generators) may be identified with the pure mapping class group of the four punctured sphere $\IP^1\setminus\{0,t,1,\infty\}$. As such it has a natural action on $\M_t$ (by pushing forward loops generating the fundamental group), and this action is the desired monodromy action.

Explicitly, upon choosing appropriate generating loops of 
$\pi_1(\IP^1\setminus\{0,t,1,\infty\})$ we see $\M_t$ may be described directly in terms of monodromy matrices:

$$\M_t \cong \left\{\ (M_1,\ldots M_4)  \ \bigl\vert\ 
M_i\in G, M_4\cdots M_1=1\right\}/G$$
which in turn is simply the quotient $G^3/G$ of three copies of 
$G$ by diagonal conjugation by $G=\SL_2(\IC)$. 
In fact this quotient has been studied classically: the ring of $G$ invariant functions on $G^3$ has $7$ generators and one relation, embedding the affine GIT quotient as a hypersurface in $\IC^7$. The particular equation for this hypersurface appears on p.366 of the book \cite{FrickeKleinI} of Fricke and Klein. The Painlev\'e VI parameters essentially specify the conjugacy classes of the four monodromies $M_i$, and serve here to fibre the six-dimensional hypersurface $G^3/G$ into a four parameter family of surfaces. Looking at the explicit equation shows they are affine {\em cubic} surfaces. In turn Iwasaki \cite{Iwas-modular} has recently pointed out that this family of cubics may be  quite simply related to Cayley's explicit family \cite{Cayley-cubic} and so contains the generic cubic surface. 

The desired action of the free group $\cF_2$ 
on the Betti spaces is given by the squares of the following ``Hurwitz" action:
$$\omega_1(M_1,M_2,M_3)= (M_2,M_2M_1M_2^{-1},M_3)$$
$$\omega_2(M_1,M_2,M_3)= (M_1,M_3,M_3M_2M_3^{-1}).$$

More explicitly if we consider simple positive loops $l_1,l_2$ in $B$ based at $\frac12$ encircling $0,1$ resp. then the monodromy of the connection on $\M_{\text{Betti}}$ around $l_i$ is given by $\omega_i^2$ (with respect to certain generators of 
$\pi_1(\IP^1\setminus\{0,\frac12,1,\infty\}$). In turn it is possible to write this action directly as an action on the ring of invariant function on $G^3$.

\section{Algebraic solutions from finite subgroups of 
$\mathbf{\SL_2(C)}$}

\subsection{What exactly is an algebraic solution?}

The obvious definition is simply an algebraic function $y(t)$ which solves \PVI\ for some value of the four parameters. Thus it will be specified by some polynomial equation
$$F(y,t)=0.$$
and a four-tuple $\theta$ of parameters.
In practice however such polynomials $F$ can be quite unwieldy and are difficult to transform under the affine Weyl symmetry group, making it difficult to see if in fact two solutions are equivalent.
This leads to our preferred definition:

{\bf Definition.\ }
An algebraic solution of \PVI\ is a compact, possibly singular, algebraic curve $\Pi$ together with two rational functions 
$y,t:\Pi\to \IP^1$

\[
\begin{CD}
\Pi @>y>> \IP^1\\
@VtVV \\
\IP^1
\end{CD}
\]

such that

$\bullet$ $t$ is a Belyi map (i.e. its branch locus is a subset of $\{0,1,\infty\}$), and

$\bullet$ $y$, when viewed as a function of $t$ away from the ramification points of $t$, solves \PVI\ for some value of the four parameters.
 
In principle it is straightforward to go between the two definitions, but in practice it is useful to look for a good model of $\Pi$ (and the model given by the closure of the zero locus of the polynomial $F$ is usually a bad choice).

\subsection{ (A) $\mapsto$ (C) }

Suppose we have a linear 
connection {\bf (A)} with finite monodromy. Its monodromy representation will be specified by a triple $(M_1,M_2,M_3)\in G^3$ generating a finite subgroup $\Gamma\subset G$ (where $G=\SL_2(\IC)$ as above). This linear connection specifies the initial value (and first derivative) 
of a solution to \PVI. This \PVI\ solution will have finite monodromy, since we know the branching of \PVI\ solutions corresponds to the $\cF_2$ action on conjugacy classes of triples in $G^3$, and the orbit through $(M_1,M_2,M_3)$ will be finite, since the action is within triples of generators of $\Gamma$.

Thus we see that finite $\cF_2$ orbits (in $G^3/G$) correspond to \PVI\ solutions with a finite number of branches, and the points of such $\cF_2$ orbits correspond to the individual branches of the \PVI\ solution. In particular the size of the orbit, the number of branches, is the degree of the map $t:\Pi\to \IP^1$.
(Indeed the $\cF_2$ action on such a finite orbit itself gives the full permutation representation of the Belyi map $t:\Pi\to \IP^1$, and in particular, by the Riemann--Hurwitz formula, determines the genus of the `Painlev\'e curve' $\Pi$.)

Said differently it is useful to define a {\em topological algebraic \PVI\ solution} (or henceforth for brevity a {\em topological solution}) to be a finite $\cF_2$ orbit in $G^3/G$. 
(The classification of such orbits is still open and is the main step in classifying all finite branching \PVI\  solutions.)
In these terms the first paragraph above points out that one obtains ``obvious" topological solutions upon taking any triple of generators of any finite subgroup of $G$.

For example (omitting discussion of how they were actually constructed) here are some solutions corresponding to certain triples of generators of the binary tetrahedral and octahedral subgroups, due to Dubrovin \cite{Dub95long} and Hitchin \cite{Hit-Octa} (in different but equivalent forms):

$$\text{Tetrahedral solution of degree three}$$
$$
y=\frac{(s-1)(s+2)}{s(s+1)},\qquad
t=\frac{(s-1)^2(s+2)}{(s+1)^2(s-2)},\qquad 
\theta=(2, 1, 1, 2)/3,
$$

$$\text{Octahedral solution of degree four}$$
$$
y=\frac{(s-1)^2}{s(s-2)},\qquad
t=\frac{(s+1)(s-1)^3}{s^3(s-2)},\qquad
\theta=(1,1,1,1)/4.
$$

In both cases $\Pi$ is a rational curve (with parameter $s$).
Although written in this compact form, one should bear in mind these formulae represent a whole (isomonodromic) family of connections {\bf (A)} as $t$ varies. An explicit elliptic solution appears in Hitchin \cite{Hit-Poncelet} and may be written as:
$$\text{Elliptic dihedral solution}\quad\theta=(1,1,1,1)/2$$
{\small
$$
y={\frac { \left( 3\,s-1 \right)  \left( {s}^{2}-4\,s-1 \right)  \left(
{s}^{2}+u \right)  \left( s \left( s+2 \right) -u \right) }{ \left( 3
\,{s}^{3}+7\,{s}^{2}+s+1 \right)  \left( {s}^{2}-u \right)  \left( s
 \left( s-2 \right) +u \right) }},\quad
t={\frac { \left( {s}^{2}+u \right) ^{2} \left( s \left( s+2 \right) -u
 \right)  \left( s \left( s-2 \right) -u \right) }{ \left( {s}^{2}-u
 \right) ^{2} \left( s \left( s+2 \right) +u \right)  \left( s \left(
s-2 \right) +u \right) }}
$$
}

where the pair $(s,u)$ lives on the elliptic curve 
$u^2=s \left( {s}^{2}+s-1 \right)$. This solution has degree $12$ and corresponds to a triple of generators of the binary dihedral group of order $20$.

It turns out (see \cite{icosa} Remark 16) that the icosahedral solutions of Dubrovin and Mazzocco \cite{DubMaz00} fit into this framework as well and correspond to (certain) triples of generators of the binary icosahedral group, although in the first instance they arose from the icosahedral reflection group as described earlier. Note that \cite{DubMaz00} Remark 0.1 describes a relation between their solutions of \PVI\ and a certain {\em folding} of Schwarz's list; this is different to the relation just mentioned---in particular problem {\bf (A)} demands an {\em extension} of Schwarz's list.

\section{Beyond Platonic $\text{\bf P}_{\text{\bf VI}}$ solutions}

My starting point in this project was simply the observation that there should be more algebraic solutions to \PVI\ than those coming from finite subgroups of $\SL_2(\IC)$. Dubrovin \cite{Dub95long} had shown how to relate  three dimensional real orthogonal reflection groups  to a certain one-parameter family 
of the full four dimensional family of \PVI\ equations (namely the family having parameters $\theta=(0,0,0,*)$) and this was used in \cite{DubMaz00} to classify  algebraic solutions having parameters in this 
one-parameter family. 
(Some aspects of \cite{DubMaz00} were subsequently extended by Mazzocco in \cite{Maz-rat} to classify rational solutions---i.e. those with only one branch, cf. also \cite{YuanLi}.)
The further observation was that if one is able to get away from the orthogonality condition here then one will relate {\em any} \PVI\  equation to a three dimensional {\em complex} reflection group.

\begin{thm}[\cite{pecr,k2p}] \label{thm: 1}
The isomonodromic deformations of type {\bf (B)} connections (on rank three vector bundles) are also controlled by the Painlev\'e VI equation, and all \PVI\ equations arise in this way.
\end{thm}

Thus a solution to \PVI\ can also be viewed as specifying an isomonodromic family of rank three connections. It turns out that the formulae to go from a \PVI\ solution $y(t)$ to such an isomonodromic family are more symmetrical than in the previous case 
(type {\bf (A)})  so we will recall them here. (For the analogous formulae for {\bf (A)} see \cite{JM81} and in Harnad's dual picture---the formula for which should be compared to that below---see \cite{Harn94} and also \cite{mazz-irreg}, which was kindly pointed out by the referee.) First the parameters: let $\la_i=\tr(B_i)$ for $i=1,2,3$ and let $\mu_i$ be the eigenvalues, in some order, of $B_1+B_2+B_3$ (which is minus the residue at infinity), so that $\sum\la_i=\sum\mu_i$.

\begin{thm}[\cite{srops}]\label{thm: 3x3}
If $y(t)$ solves Painlev\'e VI with parameters $\theta$ where
\begin{equation}\notag
\theta_1= \lambda_1-\mu_1,\quad
\theta_2= \lambda_2-\mu_1,\quad
\theta_3= \lambda_3-\mu_1,\quad
\theta_4=\mu_3-\mu_2,
\end{equation}
and we define $x(t)$ via
$$x=\frac{1}{2}\left(
\frac{(t-1)y'-\theta_1}{y}+\frac{y'-1-\theta_2}{y-t}-\frac{\, t\,y'+\theta_3}{y-1}
\right)
$$
then the family of logarithmic connections {\bf (B)}
will be isomonodromic as $t$ varies, where

$$  
B_1=\left(\begin{matrix}
\lambda_1& b_{12} & b_{13} \\
0&0&0\\
0&0&0
\end{matrix}\right),
\quad
B_2=\left(\begin{matrix}
0&0&0 \\
b_{21}&\lambda_2&b_{23} \\
0&0&0
\end{matrix}\right),\quad
B_3=\left(\begin{matrix}
0&0&0 \\
0&0&0 \\
b_{31}&b_{32}&\lambda_3
\end{matrix}\right)
$$

\begin{align*}
b_{12}&=
{ \lambda_1}\,\, - \mu_3y +(\mu_1 - xy)(y-1),
\ &
b_{32}&=(\mu_2-\lambda_2-b_{12})/t,
\\
b_{13}&=
  { \lambda_1}t-\mu_3y +(\mu_1 - xy)(y-t),
\ &
b_{23}&=(\mu_2-\lambda_3)t-b_{13},
\\
b_{21}&=
\lambda_2+\frac{\mu_3(y-t)-\mu_1(y-1)+x(y-t)(y-1)}{t-1},
\ &
b_{31}&=(\mu_2-\lambda_1-b_{21})/t.\\
\end{align*}
\end{thm}

The implication of this for algebraic solutions should now be clear: the monodromy of a \PVI\ solution is also described by an action of the free group $\cF_2$ on (conjugacy classes of) triples of three dimensional complex reflections $(r_1,r_2,r_3)$ (with the same formula as before, just replace $M_i$ by $r_i$). Thus in this context the ``obvious" topological solutions (i.e finite $\cF_2$ orbits) come from taking a triple of generating reflections of a finite complex reflection group in 
$\GL_3(\IC)$. Such finite complex reflection groups were classified in 1954 by Shephard and Todd \cite{Shep-Todd} and apart from the familiar real orthogonal reflection groups there is an infinite family plus four exceptional complex groups, the Klein reflection group (of order $336$, a two-fold cover of Klein's simple group isomorphic to $\PSL_2(\IF_7)\hookrightarrow\PGL_3(\IC)$), two Hessian groups and the Valentiner group (of order 2160, a six fold cover of $A_6\hookrightarrow\PGL_3(\IC)$).

The infinite family of groups and the two Hessian groups do not seem to lead to interesting new solutions, but by computing the $\cF_2$ orbits (determining the topology of $\Pi$) it is easy to see that the Klein group yields a genus zero degree $7$ solution and the Valentiner group has three inequivalent triples of generating reflections, each leading to genus one solutions with degrees $15,15,24$ respectively. 
These are new solutions, previously undetected.
(The $24$ appearing here led to a certain amount of trepidation, given that the $10$ page elliptic solution of \cite{DubMaz00} had degree $18$.)

\subsection{Construction} 

Of course finding the topological solution is not the same as finding an explicit isomonodromic family of connections; one needs to solve a family of Riemann--Hilbert problems inverting the transcendental Riemann--Hilbert map for each value of $t$.
(Indeed the author's original plan was to just prove the existence of new interesting solutions, in \cite{pecr}, but a certain stubbornness, and some inspiration from reading about Klein's work finding explicit $3\times 3$ matrices generating his simple group, convinced us to go further.)

The two main steps in the method we finally got to work are as follows. 
(This is a generalisation of the method used by Dubrovin--Mazzocco \cite{DubMaz00}.)

{\bf 1)  Jimbo's asymptotic formulae.}\ \ 
In \cite{Jimbo82} M. Jimbo found an exact formula for the leading asymptotics at $t=0$ of the branch of the \PVI\  solution $y(t)$ corresponding to any sufficiently generic linear monodromy representation $(M_1,M_2,M_3)$.
(This formula was obtained by considering the degeneration of the isomonodromic family of connections {\bf (A)} as $t\to 0$; in the limit the four punctured sphere degenerates into a stable curve with two components, each with three marked points. The connections {\bf (A)} degenerate into hypergeometric systems on each component, with known monodromy. Since these are rigid it is easy to solve their Riemann--Hilbert problems explicitly and this gives the leading asymptotics of the isomonodromic family and thus of the \PVI\ solution.)

\ 

\begin{figure}[h]
	\centering
	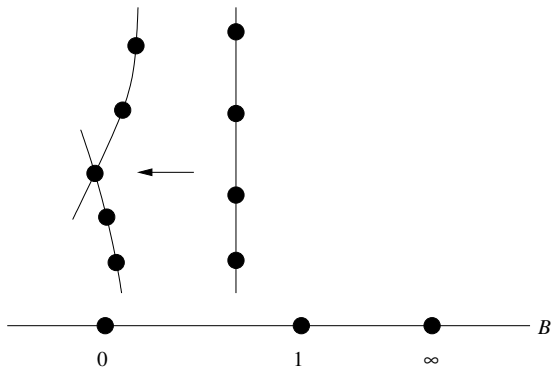
	\caption{Degeneration to two hypergeometric systems}
\end{figure}

This is useful for us because, as Jimbo mentions, one may substitute the leading asymptotics back into the \PVI\ equation to get arbitrarily many terms of the precise asymptotic expansion of the solution at $0$. If the solution is algebraic, then this is its Puiseux expansion, a sufficient number of terms of which will determine the entire solution.

It turns out there was a typo in \cite{Jimbo82}, which meant the entire method did not work (indeed the fact it did not work led to the questioning of Jimbo's formula and hence the correction \cite{k2p}). (Note the special parameters of \cite{DubMaz00} are not covered by Jimbo's result; rather \cite{DubMaz00} adapted the argument of \cite{Jimbo82} to their case.) 

{\bf 2) Relating (A) and (B).}\ \ 
Since Jimbo's formula requires a monodromy representation of a connection of type {\bf (A)}, and we are starting with a triple of $3\times 3$ complex reflections (the monodromy representation of a connection of type {\bf (B)}), the second step is that we need to see how to go between these two pictures (on both the De Rham and Betti sides of the Riemann--Hilbert correspondence).
This will be described in the following subsection.

\subsection{Relating connections {\bf (A)} and {\bf (B)}}

We wish to sketch how to convert a connection {\bf (B)} on a rank three vector bundle into a connection of the form {\bf (A)} on a rank two bundle. On the other side of the Riemann--Hilbert correspondence this amounts to an $\cF_2$-equivariant map from triples of complex reflections to triples of elements of $G=\SL_2(\IC)$ (as in \cite{k2p} Section 2).

Of course the monodromy groups change in a highly nontrivial way under this procedure. For example the Klein reflection group becomes the triangle group $\Delta_{237}\subset G$, which is an infinite group, and the Valentiner group becomes the binary icosahedral group (leading to an unexpected relation between $A_6$ and $A_5$).

After this procedure was put on the arxiv (\cite{k2p}) we learnt \cite{Dett-Reit-p6} that it is essentially a case of N. Katz's middle convolution functor \cite{Katz-rls}, although our construction using the complex analytic Fourier--Laplace transform is different from that of Katz (using l-adic methods) and from the work of Dettweiler--Reiter \cite{DettReit-katz}. 

The basic picture which emerges is as follows (see Diagram 1), and ought to be better known. It was obtained essentially by a careful reading of 
the article \cite{BJL81} of Balser, Jurkat and Lutz, although the basic idea of relating irregular and Fuchsian systems by the Laplace transform dates back to Birkhoff and Poincar\'e.
(Dubrovin \cite{Dub95long, DubPT} used an orthogonal analogue in relation to Frobenius manifolds, also using \cite{BJL81}. Moreover the top triangle is essentially a case of `Harnad duality' \cite{Harn94} so for $n=3$ we knew we would obtain all \PVI\ equations.)

\begin{equation*}
 \setlength{\unitlength}{60pt}
 \begin{picture}(4.2,4.1)(0,0)
 \put(2,-0.5){\object{Diagram 1}}
 \put(2,3.9){\object{$A\in \gl_n(\IC)$}}
 \put(0,3){\object{$d-\sum\frac{A_i}{z-a_i}dz$}}
 \put(4,3){\object{$d-\left(\frac{A^0}{z^2}+\frac{A}{z}\right)dz$}}
 \put(0,1.1){\object{$r_1,\ldots,r_n$}}
 \put(0,0.9){\object{\footnotesize $r_i=1+e_i\otimes\al_i$}}
 \put(4,1){\object{$\qquad\qquad\qquad(u_-,u_+,h)\in U_-\times U_+\times H$}}
 \put(2,0.1){\object{$\GL_n(\IC)$}}
 \put(2.5,3.75){\vector(2,-1){1}}
 \put(1.5,3.75){\vector(-2,-1){1}}
     \put(1.1,3.3){\selabel{\footnotesize $A_i=E_iA$}}
 \put(0,2.75){\vector(0,-1){1.4}}
    \put(0.7,1.9){\elabel{\small Riemann--Hilbert}}
 \put(4,2.75){\vector(0,-1){1.4}}
    \put(3.7,1.9){\elabel{\small Stokes}}
 \put(0.8,3){\vector(1,0){2.4}}
 \put(3.2,3){\vector(-1,0){2.4}}
	\put(1.8,2.85){\slabel{\small Fourier--Laplace}}
 \put(0.8,1){\vector(1,0){2.4}}
 \put(3.2,1){\vector(-1,0){2.4}}
	\put(1.8,1.15){\slabel{\small Killing--Coxeter}}
	\put(1.8,0.85){\slabel{\footnotesize $(\al_i(e_j))=hu_+-u_-$}}
 \put(0.5,0.75){\vector(2,-1){1}}
	\put(0.6,0.5){\swlabel{\small $r_n\cdots r_1$}}
 \put(3.5,0.75){\vector(-2,-1){1}}
    \put(3.2,0.5){\swlabel{\small $u_-^{-1}hu_+$}}
 \end{picture}
\end{equation*}

\ 

\

The basic idea is to describe a transcendental map from from $\gl_n(\IC)$ to $\GL_n(\IC)$ in two different ways (the two paths down the left and the right from the top to the bottom of the diagram).

Choose $n$ distinct complex numbers $a_1,\ldots,a_n$ and define $A^0=\diag(a_1,\ldots,a_n)$. Roughly speaking (on a dense open patch) the left-hand column arises by defining $A_i=E_iA$ (setting to zero all but the $i$th row of $A$) and constructing the logarithmic connection $d-\sum\frac{A_i}{z-a_i}dz$ having rank one residues at each $a_i$. Then taking the monodromy of this yields $n$ complex reflections $r_i$ (and if bases of solutions are chosen carefully one can naturally define vectors $e_i$ and one-forms $\al_i$ such that $r_i=1+e_i\otimes \al_i$ and that the $e_i$ form a basis). Then  the map to $\GL_n(\IC)$ is given by taken the product of $r_n\ldots r_1$ of these reflections, written in the $e_i$ basis.

Now the key algebraic fact, which dates back at least to Killing \cite{killing2} (see Coleman
\cite{Coleman}), is that any such product of complex reflection lies in the big cell of $\GL_n(\IC)$ and so may be factored as the product of a lower triangular and an upper triangular matrix. We write this product as $u_-^{-1}hu_+$ with $u_\pm\in U_\pm$ the unipotent triangular subgroups, and $h\in H$ diagonal:
\beq \label{eq: KC matrix}
r_n\cdots r_2r_1 = u_-^{-1}hu_+ .\eeq
Further, although this relation between the reflections and $u_\pm$ looks to be highly nonlinear, one can relate them in an almost linear fashion:
the matrix $hu_+-u_-$ is the matrix with entries $\al_i(e_j)$.

On the other hand it turns out that the same map can be defined by taking the {\em Stokes data} of the irregular connection 
$d-\left(\frac{A^0}{z^2}+\frac{A}{z}\right)dz$ (indeed the map on the right-hand side generalises \cite{bafi} to any complex reductive group $G$ in place of $\GL_n(\IC)$, but only for $\GL_n(\IC)$ is the alternative ``logarithmic" viewpoint available).
Thus $u_\pm$ are also the two Stokes matrices of this irregular connection (the natural analogue of monodromy data for such connections);
the exact definition is not important here. (The element $h$ is the so-called formal monodromy, explicitly it is simply $\exp(2\pi i \Lambda)$ where $\Lambda$ is the diagonal part of $A$.)
The two connections are related  (see \cite{BJL81}) by the Fourier--Laplace transform: this is more than just formal, and by relating bases of solutions on both sides the stated relation between the Stokes and monodromy data is obtained. (In both cases the resulting element of $\GL_n(\IC)$ is the monodromy around $z=\infty$ in a suitable basis.) In summary we see that the ``Betti" incarnation of the Fourier--Laplace transform is the relation of Killing--Coxeter.

Now to apply this in the current context we consider the effect of adding a scalar $\la$ to $A\in\gl_n(\IC)$. On the right-hand side this corresponds to tensoring the irregular connection by the logarithmic connection $d-\la dz/z$ on the trivial line bundle, and it follows \cite{BJL81} that the Stokes data is changed only by scaling $h$ by $s:=\exp(2\pi i \la)$, fixing $u_\pm$.
On the logarithmic side this corresponds to a nontrivial convolution operation, changing the monodromy representation in a nontrivial way. Of course using the Killing--Coxeter identity we now see precisely how the complex reflections vary. (It is perhaps worth noting that this scalar shift is essentially the inverse of the spectral parameter introduced by Killing
\cite{killing2} p.20, appearing in the characteristic polynomial of the Killing--Coxeter matrix 
\eqref{eq: KC matrix}: $\det(u_-^{-1}shu_+-1)=\det(shu_+-u_-)$.)

If we set $n=3$ then  the logarithmic connections appearing are 
of the form {\bf (B)}, upon taking $a_1,a_2,a_3=0,t,1$.
Then we may choose the scalar shift such that the resulting element of 
$\GL_3(\IC)$ has $1$ as an eigenvalue. This implies that the connections are reducible and we can take the irreducible rank two sub- or quotient connection. Projecting to $\lsl_2$ gives the desired connection of type {\bf (A)} (see \cite{k2p}).
(Note that there is a choice involved here, of which eigenvalue to shift to $1$.) 

\subsection{New solutions}

Thus in summary the procedure now is as follows: take a triple of generating reflections of a finite complex reflection group in $\GL_3(\IC)$. Push it down to the $2\times 2$ framework using the scalar shift to obtain a triple $(M_1,M_2,M_3)$ of elements of $\SL_2(\IC)$ in an isomorphic $\cF_2$ orbit. Apply Jimbo's formula to get the leading asymptotics of the corresponding \PVI\ solutions at $t=0$ on each branch (i.e. for each triple in the $\cF_2$ orbit). (Converting the values which arise into exact algebraic numbers.)
Substitute these leading terms back into \PVI\ to obtain arbitrarily many terms of the Puiseux expansion at $0$ of each solution branch. Use these expansions to determine the polynomial $F(y,t)$ defining the solution (assuming it is algebraic).
Find a parametrisation of the resulting algebraic curve (for example using M. van Hoeij's wonderful Maple algebraic curves package).

For example for the Klein complex reflection group of order $336$ this works perfectly (\cite{k2p}) and the resulting solution is:
$$\text{Klein solution  }\quad \theta=(2,2,2,4)/7$$
$$
y=-{\frac { \left( 5\,{s}^{2}-8\,s+5 \right)  \left( 7\,{s}^{2}-7\,s+4\right) }
{ s \left( s-2 \right)  \left( s+1 \right) \left( 2\,s-1 \right) 
\left( 4\,{s}^{2}-7\,s+7 \right)  }}\qquad
t={\frac { \left( 7\,{s}^{2}-7\,s+4 \right) ^{2}}{{s}^{3} \left( 4\,{s}^
{2}-7\,s+7 \right) ^{2}}}$$

which has $7$ branches. One may of course now substitute this back into the formula of Theorem \ref{thm: 3x3} (with $\la=(1,1,1)/2, \  \mu=(3,5,13)/14)$  to obtain an explicit family of logarithmic connections having monodromy equal to the Klein reflection group generated by reflections
(see \cite{srops} \S3).

When converted to connections of type {\bf (A)} these ``Klein connections" have infinite (projective) monodromy group equal to the triangle group $\Delta_{237}$ (cf. \cite{octa} Appendix B). 
On the other hand it turns out \cite{icosa} that for the Valentiner connections, even though they are much trickier to construct directly, we can still compute immediately that they  become connections of type {\bf (A)} with binary icosahedral monodromy. They are also inequivalent  to those appearing in the work of Dubrovin--Mazzocco related to the real orthogonal icosahedral reflection group (which lead to unipotently generated  monodromy with one choice of the scalar shift, but finite binary icosahedral monodromy with a different choice, cf. \cite{icosa} Remark 16).

Thus it seemed like a good idea to examine  precisely what \PVI\ solutions arise upon taking arbitrary triples of generators $(M_1,M_2,M_3)$ of the binary icosahedral group. Thus we looked at all triples of generators and quotiented by the relation coming from the affine $F_4$ symmetries of \PVI.
The resulting table has $52$ rows (which is quite small considering there are $26688$ conjugacy classes of generating triples). The first $10$ rows correspond to the $10$ icosahedral rows of Schwarz's list and thus the projective monodromy around one of the four punctures is the identity (these correspond to the \PVI\ solution $y=t$). The remaining rows are as in Table \ref{table: icosalist} (this is abridged from \cite{icosa}). (Note that the right notion of equivalence in the linear nonrigid problem {\bf (A)} seems to be the `geometric equivalence' of \cite{icosa} section 4----however this coincides with equivalence under the affine $F_4$ Weyl group, in this case.)

\begin{table}[h]
\begin{center}
\begin{tabular}{|c|c|c|c|c|||c|c|c|c|c| }
\hline
  & \text{Degree}  &  \text{Genus} & \text{Walls} & \text{Type}  & &\text{Degree}  &  \text{Genus} & \text{Walls} & \text{Type} \\ \hline
11 & 2 & 0 & 2 & $ b^2\,c^2$ & 32 & 10 & 0 & 3 & $ d^4$ \\ \hline
12 & 2 & 0 & 2 & $ b^2\,d^2$ &33 & 12 & 0 & 0 & $ a\,b\,c\,d$ \\ \hline
13 & 2 & 0 & 2 & $ c^2\,d^2$ &34 & 12 & 1 & 1 & $ a\,b\,c^2$ \\ \hline
14 & 3 & 0 & 1 & $ b\,c^2\,d$&35 & 12 & 1 & 1 & $ a\,b\,d^2$ \\ \hline
15 & 3 & 0 & 1 & $ b\,c\,d^2$&36 & 12 & 1 & 1 & $ b^2\,c\,d$ \\ \hline
16 & 4 & 0 & 2 & $ a\,c^3$ &37 & 15 & 1 & 2 & $ b^3\,c$ \\ \hline
17 & 4 & 0 & 2 & $ a\,d^3$ &38 & 15 & 1 & 2 & $ b^3\,d$ \\ \hline
18 & 4 & 0 & 2 & $ c^3\,d$ &39 & 15 & 1 & 2 & $ b^2\,c^2$ \\ \hline
19 & 4 & 0 & 2 & $ c\,d^3$ &40 & 15 & 1 & 2 & $ b^2\,d^2$ \\ \hline
20 & 5 & 0 & 1 & $ b^2\,c\,d$&41 & 18 & 1 & 3 & $ b^4$  \\ \hline
21 & 5 & 0 & 2 & $ c^2\,d^2$ &42 & 20 & 1 & 1 & $ a\,b^2\,c$ \\ \hline
22 & 6 & 0 & 1 & $ b\,c^2\,d$&43 & 20 & 1 & 1 & $ a\,b^2\,d$ \\ \hline
23 & 6 & 0 & 1 & $ b\,c\,d^2$&44 & 20 & 1 & 3 & $ a^2\,c^2$ \\ \hline
24 & 8 & 0 & 1 & $ a\,c^2\,d$&45 & 20 & 1 & 3 & $ a^2\,d^2$ \\ \hline
25 & 8 & 0 & 1 & $ a\,c\,d^2$&46 & 24 & 1 & 2 & $ a\,b^3$ \\ \hline
26 & 9 & 1 & 2 & $ b\,c^3$ &47 & 30 & 2 & 2 & $ a^2\,b\,c$ \\ \hline
27 & 9 & 1 & 2 & $ b\,d^3$ &48 & 30 & 2 & 2 & $ a^2\,b\,d$ \\ \hline
28 & 10 & 0 & 2 & $ a^2\,c\,d$&49 & 36 & 3 & 3 & $ a^2\,b^2$ \\ \hline
29 & 10 & 0 & 2 & $ b^3\,c$ &50 & 40 & 3 & 3 & $ a^3\,c$ \\ \hline
30 & 10 & 0 & 2 & $ b^3\,d$ &51 & 40 & 3 & 3 & $ a^3\,d$ \\ \hline
31 & 10 & 0 & 3 & $ c^4$ &52 & 72 & 7 & 3 & $ a^3\,b$\\\hline
\end{tabular}
\caption{Icosahedral solutions $11-52$}\label{table: icosalist}
\end{center}\end{table}

Thus there are lots of other icosahedral solutions the largest having genus $7$ and $72$ branches. (The column ``Type'' indicates the set of conjugacy classes of local monodromy of the corresponding connections of type {\bf (A)}, as we marked on Schwarz's list. The column ``Walls'' indicates the number of reflection hyperplanes for the affine $F_4$ Weyl group that the solution's parameters $\theta$ lie on.) 
A few of these solutions had appeared before:
those with degree less than $5$
are simple deformations of previous solutions, solutions 21 and 26 are in Kitaev \cite{Kitaev-dessins} and the Dubrovin--Mazzocco icosahedral solutions are equivalent to those on rows $31,32,41$. 
On the other hand the Valentiner solutions are quite far down the list on rows $37,38$ and $46$.

The above method of constructing solutions using Jimbo's asymptotic formula applies only to sufficiently generic monodromy representations but it turns out that most of the rows of this table have some representative (in their affine $F_4$ orbit) to which Jimbo's formula maybe applied (on every branch).
Thus we could start working down the list constructing new solutions. An initial goal was to get to solution $33$: this solution purports to be on {\em none} of the reflection hyperplanes and the folklore was that all explicit solutions to Painlev\'e equations must lie on some reflection hyperplane. The folklore was wrong:
\begin{gather}
\text{``Generic'' solution,\quad  row $33$,\quad  $\theta=(2/5, 1/2,1/3, 4/5)$}\notag\\
y=
-\frac{9 s (s^2+1) (3 s-4) (15 s^4-5 s^3+3 s^2-3 s+2)}
{(2 s-1)^2(9 s^2+4)(9 s^2+3 s+10)}\qquad
t=\frac{27 s^5 (s^2+1)^2 (3 s-4)^3}{4(2 s-1)^3(9 s^2+4)^2}.\notag
\end{gather}
So far this looks to be  the only example of a `classical' solution of any of the Painlev\'e equations that does not lie on a reflection hyperplane (of the full symmetry group). 
Apart from being in the interior of a Weyl alcove this solution is generic in another sense: a randomly chosen triple of generators of the binary icosahedral group is most likely to lead to it (more of the $26688$ triples of generators correspond to this row than to any other). Notice also that this solution has type $abcd$; there is one local monodromy in each of the 
four nontrivial conjugacy classes of $A_5$.

At this stage we were approaching solution $41$ which we knew took $10$ pages to write down. So we stopped and looked around to see if there were other interesting (even just topological) solutions. (The tetra/octahedral cases could all now be fully dealt with \cite{octa}.)

\section{Pullbacks}

In his book \cite{Klein-ico} on the icosahedron Klein showed that 
all second order Fuchsian differential equations with finite monodromy are (essentially) pullbacks of a hypergeometric equation
along a rational map $f$:

\begin{equation*}
 \setlength{\unitlength}{60pt}
 \begin{picture}(4.2,1.2)(0,0)
 \put(0,1){\object{$\IP^1\setminus\{\text{$k$ points}\}$}}
 \put(3.7,1){\object{$\IP^1$}}
 \put(2,0){\object{$\IP^1\setminus\{0,1,\infty\}$}}
 \put(0.8,1){\line(1,0){0.2}}
 \put(1.2,1){\line(1,0){0.2}}
 \put(1.6,1){\line(1,0){0.2}}
 \put(2.0,1){\line(1,0){0.2}}
 \put(2.4,1){\line(1,0){0.2}}
 \put(2.8,1){\line(1,0){0.2}}
 \put(3.2,1){\vector(1,0){0.2}}
	\put(1.8,1.15){\slabel{\small Schwarz map}}
	\put(1.8,0.85){\slabel{\footnotesize $y_1/y_2$}}
 \put(0.5,0.75){\vector(2,-1){1}}
	\put(0.6,0.5){\swlabel{$f$}}
 \put(3.5,0.75){\vector(-2,-1){1}}
    \put(3.4,0.5){\swlabel{\small $/\Gamma$ (invariants)}}
 \end{picture}
\end{equation*}

In particular ($k=3$) all the icosahedral entries on Schwarz's list, may be obtained by pulling back the``$235$" hypergeometric equation (on row VI of Schwarz's list).

In our context, an isomonodromic family of connections of type {\bf (A)} amounts to a family of Fuchsian equations with $5$ singularities (at $0,t,1,\infty$ plus an apparent singularity at another point $y$).\footnote{this is the same $y$ appearing in \PVI---i.e. the function $y$ on the space of connections {\bf (A)} is the position of the apparent singularity that appears when the connection is converted into a Fuchsian equation \cite{RFuchs}.} 
Klein's theorem says each element of this family arises as the pullback  
of the $235$ hypergeometric equation along a rational map, so the family corresponds to a family of rational maps.

Thus finding a \PVI\ solution corresponding to a family of connections {\bf (A)} with finite monodromy amounts to giving a certain family of rational maps $f:\IP^1\to\IP^1$.
To construct such \PVI\ solutions one may try to find such families of rational maps, such that each map pulls back a hypergeometric equation to an equation with the right number of singular points---or to one that can be put in this form after using elementary transformations to remove extraneous apparent singularities. (This is not straightforward; for example given a finite monodromy representation of a connection {\bf (A)} it is not immediate even what degree such a map $f$ will have.)

An important further observation (due to C. Doran \cite{Chuck1} and A. Kitaev \cite{Kit-sfit6}) is that any such family of rational maps will lead to algebraic solutions of Painlev\'e VI regardless of whether or not the hypergeometric equation being pulled back has finite monodromy (provided the equation upstairs has the right number of poles); the algebraicity follows from that of the family of rational maps.

Kitaev and Andreev \cite{Kit-sfit6, And-Kit-CMP, Kitaev-dessins} have used this to construct some \PVI\ solutions, essentially by starting  to enumerate all such rational maps (this leads to a few new solutions, but most in fact turn out to be equivalent to each other or to ones previously constructed---see the summary at the end of this article).

On the other hand, Doran had the idea that interesting \PVI\ solutions should come from hypergeometric equations with interesting 
monodromy groups. Thus (amongst other things) \cite{Chuck1} studied the possible hypergeometric equations with monodromy a {\em hyperbolic arithmetic triangle group} which may be pulled back to yield \PVI\ solutions. Indeed in \cite{Chuck1} Corollary 4.6, Doran lists such possible triangle groups and the degrees and ramification indices of the corresponding rational maps $f$, although no new solutions were actually constructed. 
We picked up on this thread in \cite{octa} Section 5: it was found that all but one entry on Doran's list corresponded to a known explicit solution (although were perhaps unknown when \cite{Chuck1} was published). The remaining entry was for a family of degree $10$ rational maps $f$ pulling back the $237$ triangle group with ramification indices (partitions of $10$):
$$[2,2,2,2,2],\ [3,3,3,1],\ [7,1,1,1]$$
over $0,1,\infty$ (where the hypergeometric system has projective monodromy of orders $2,3,7$ resp.), as well as minimal ramification $[1^8,2]$ over another variable point.
As explained in \cite{octa} one can get from here to a topological \PVI\  solution by drawing a picture: we wish to find such a rational map $f$ topologically---i.e. describe the topology of a branched cover $f:\IP^1\to\IP^1$  with this ramification data. This may be done by playing ``join the dots" (completely in the spirit of Grothendieck's Dessins d'Enfants) and yields a covering diagram as required. One diagram so obtained is shown in figure \ref{fig: dessin}. (Note that the idea of drawing pictures such
     as figure \ref{fig: dessin} first appeared in Kitaev \cite{Kitaev-dessins}.)

\begin{figure}[h] 
	\centering
	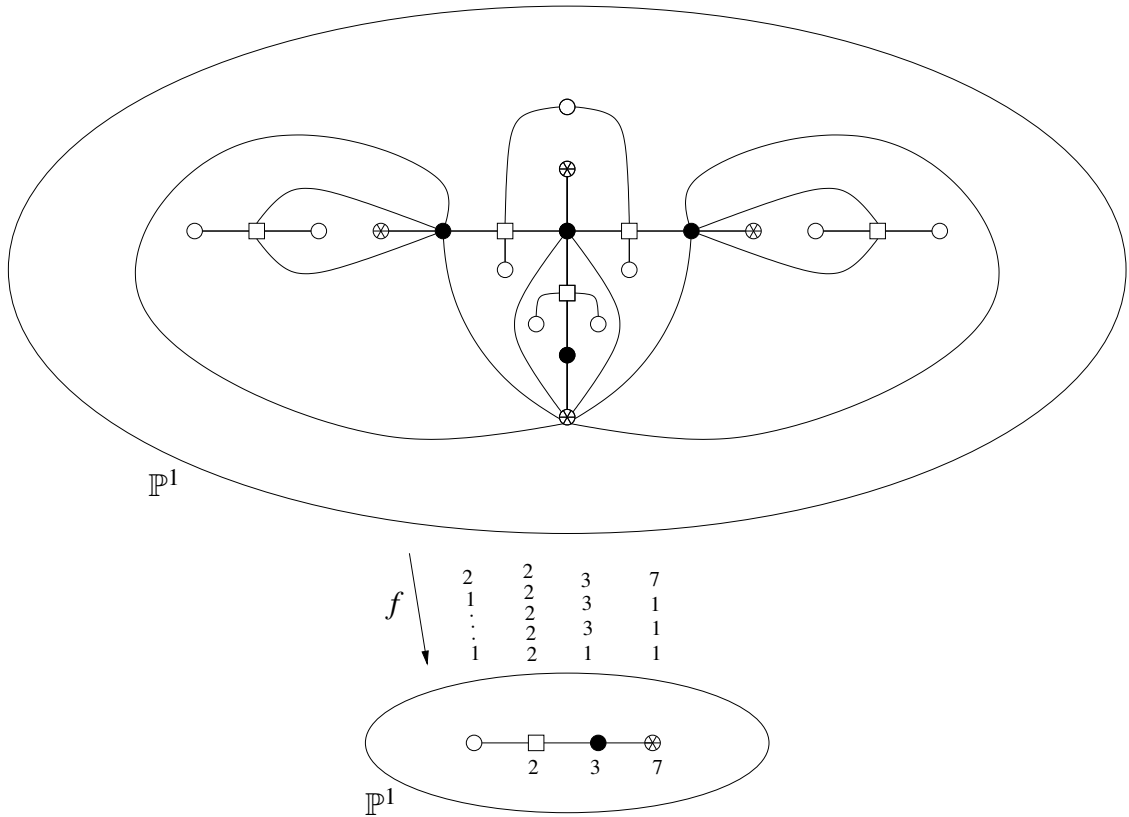
	\caption{$237$ degree $10$ rational map $f$}\label{fig: dessin}
\end{figure}

The upper copy of $\IP^1$ is thus divided into $10$ connected components and $f$ maps each component isomorphically onto the complement of the interval drawn on the lower $\IP^1$ (the lines and the vertices upstairs are the preimages of the lines and vertices downstairs). In particular the diagram shows how loops upstairs map to words in the generators of the fundamental group $\pi_1(\IP^1\setminus\{0,1,\infty\})$ downstairs.
In this way we can compute by hand the monodromy of the equation upstairs obtained by pulling back a hypergeometric equation with monodromy $\Delta_{237}$. This yields the triple: 
$$
M_1=caca^{-1}c^{-1},\quad 
M_2=c,\quad 
M_3=c^{-1}a^{-1}cac$$
(where $a,b,c$ are lifts to $\SL_2(\IC)$ of standard generators of $\Delta_{237}$ with $cba=1$), which 
we know a priori lives in  a finite $\cF_2$ orbit.
One finds immediately that the orbit through the conjugacy class of this triple has size $18$ and constitutes a genus one, degree $18$ topological \PVI\ solution.

Now it turns out that Jimbo's formula may be applied to every branch of this solution, and proceeding as before we obtain the solution explicitly:
\begin{gather}
\text{Elliptic $237$ solution}\qquad
\theta =(2/7, 2/7, 2/7, 1/3)\notag\\
y=
\frac{1}{2}-{\frac { \left( 3\,{s}^{8}-2\,{s}^{7}-4\,{s}^{6}-204\,{s}^{5}
-536\,{s}^{4}-1738\,{s}^{3}-5064\,{s}^{2}-4808\,s-3199 \right) u}
{4\, \left( {s}^{6}+196\,{s}^{3}+189\,{s}^{2}+756\,s+154 \right)  \left( {
s}^{2}+s+7 \right)  \left( s+1 \right) }}\notag \\\notag \\
t=
\frac{1}{2}-{\frac { \left( {s}^{9}-84\,{s}^{6}-378\,{s}^{5}
-1512\,{s}^{4}-5208\,{s}^{3}-7236\,{s}^{2}-8127\,s-784 \right) u}
{432\,s \left( s+1 \right) ^{2} \left( {s}^{2}+s+7 \right) ^{2}}}\notag
\end{gather}
where 
$u^2=s\,(s^2+s+7).$ (This solution, or rather an inequivalent `Galois conjugate' of it, has also been obtained independently by A. Kitaev \cite{kitaev-rmks} p.219 by directly computing such a family of rational maps---apparently also influenced by Doran's list.)

\section{Final steps}

\subsection{Up to degree $\mathbf{24}$}

We now have an example of a degree $18$ elliptic solution to \PVI\ with a quite simple form. This leads immediately to the suspicion that the $10$ page Dubrovin--Mazzocco solution is just written at a bad value of the parameters. Indeed using the method we have been `tweaking' while working down the icosahedral table enables us to guess good a priori choices of the parameters $\theta$ within the affine $F_4$ equivalence class 
for row 41 in Table \ref{table: icosalist}   (i.e. so that the expression for the polynomial $F$ will be `small'). Choosing such parameters and constructing the solution from scratch at those parameters yields:
\begin{thm}[\cite{icosa}]
The Dubrovin--Mazzocco icosahedral solution is equivalent to the solution 
$$
\ \qquad y=
\frac{1}{2}-
{\frac {8\,{s}^{7}-28\,{s}^{6}+75\,{s}^{5}+31\,{s}^{4}-269\,
{s}^{3}+318\,{s}^{2}-166\,s+56}
{18\,u \left( s-1 \right)  \left( 3\,{s}^{3}-4\,{s}^{2}+4\,s+2 \right) }}$$
$$$$ %
$$t=
\frac{1}{2}+
{\frac { \left( s+1 \right)  
\left( 
32\,({s}^{8}+1)-320\,({s}^{7}+s)+1112\,({s}^{6}+s^2)
-2420\,({s}^{5}+s^3)+3167\,{s}^{4}
\right) }
{54\,{u}^{3}\,s\, \left( s-1 \right) }}
$$%
\noindent
on the elliptic curve
$$\ \qquad u^2=s\,(8\,s^2-11\,s+8)$$
\noindent
with $\theta=(1,1,1,1)/3$. In particular this elliptic curve is birational to that defined by the ten page polynomial.
\end{thm}
Substituting this into the formula of Theorem \ref{thm: 3x3} with
$\lambda=(1,1,1)/2,\mu=(1,3,5)/6$ now gives the third (and trickiest) family of connections of type {\bf (B)} with monodromy the icosahedral reflection group.

This can be pushed further with more tweaking to get up to degree $24$ (row 46 in Table \ref{table: icosalist})  i.e. to obtain the largest Valentiner solution \cite{icosa} 
(the main further tricks used are described in \cite{octa} Appendix C).
In particular this finishes the construction of all {\em elliptic} icosahedral solutions. Intriguingly, one finds that  
the resulting elliptic icosahedral Painlev\'e curves $\Pi$ become singular only 
on reduction modulo the primes $2$, $3$ and $5$ (except for rows $44,45$---we will see another reason in the following subsection that these are abnormal). Similarly the elliptic Painlev\'e curve related to the $237$ triangle group becomes singular only on reduction modulo $2$, $3$ and $7$.

\subsection{Quadratic/Landen/Folding transformations}

Now the happy fact is that the remaining icosahedral solutions may be obtained from earlier solutions by a trick, first introduced  in the context of \PVI\   by Kitaev \cite{Kitaev-quad-p6} and a simpler equivalent form  was found by Ramani et al \cite{RGT-quad}. Manin \cite{Manin-P6} refers to some equivalent transformations as Landen transformations (Landen has clear precedence since the original Landen transformations were rediscovered by Gauss!). Tsuda et al \cite{TOS-folding} call them folding transformations.

In any case the basic idea is simple: if one has a connection {\bf(A)} with two local projective monodromies of order two (say at $0,\infty$) then one can pull it back along the map $z\mapsto z^2$ and obtain a connection with only apparent singularities at $0,\infty$ (which can be removed) and four genuine singularities. This can be normalised into the form {\bf (A)}, and the key point is that this works in families and maps isomonodromic deformations of the original connections to isomonodromic deformations of the resulting connections---i.e. it transforms certain solutions of  \PVI\ into different, generally inequivalent, solutions. Of course this is not a genuine symmetry of \PVI\ since special parameters are required, but it is precisely what is needed to construct the remaining solutions. 

Indeed observe that each of the rows of the icosahedral table with degree greater than $24$ have type $a^2\xi\eta$ for some $\xi,\eta\in\{a,b,c,d\}$---i.e. they have two projective monodromies of order two. Pulling back along the squaring map will transform the corresponding connections into connections of type $\xi^2\eta^2$. It turns out (in this icosahedral case) the corresponding \PVI\ solutions have half the degree, and we obtain an algebraic relation between the solutions.
This program is carried out in \cite{ipc} and the remaining icosahedral solutions are obtained. See also Kitaev--Vid\={u}nas \cite{kit-vid-quad}. (Notice also that the elliptic solutions on rows $44,45$ are related in this way to earlier solutions.) For example in \cite{ipc} we found a relatively simple  explicit equation for the genus $7$ algebraic curve naturally attached to the icosahedron, on which the  largest (degree $72$) icosahedral solution is defined: it may be modelled as the plane octic with affine equation

$$
9\,(p^6\,q^2+p^2\,q^6)+
18\,p^4\,q^4+$$$$
4\,(p^6+q^6)+
26\,(p^4\,q^2+p^2\,q^4)+
8\,(p^4+q^4)+
57\,p^2\,q^2+$$$$
20\,(p^2+q^2)+
16=0.
$$
$$\text{The genus seven icosahedral Painlev\'e curve}$$

\section{Conclusion}

Thus in conclusion we have filled in a number of rows of what could be called the {\em nonlinear Schwarz's list}. Whether or not there will be other rows remains to be seen. So far this list of known algebraic solutions to \PVI\  takes the following shape (we will use the letters $d$ and $g$ to denote the degree and genus of solutions, and consider solutions up to equivalence under Okamoto's affine $F_4$ symmetry group. Some non-trivial work has been done to establish which of the published solutions are equivalent to each other and which were genuinely new):

First there are the {\bf rational solutions} ($d=1$), studied by Mazzocco \cite{Maz-rat} and Yuan--Li \cite{YuanLi}, which fit into the set of Riccati solutions classified by Watanabe \cite{watanabePVI}. (Beware that `rational' 
here means the solution is a rational function of $t$, which implies, but is by no means equivalent to, having a rational parameterisation.)

Then there are {\bf three continuous families} of solutions $g=0, d=2,3,4$. 
The degree two family is $y=\sqrt{t}$ which, as one may readily verify, solves \PVI\ for a family of possible parameter values. Similarly the degree $3$ tetrahedral solution, and the degree $4$ octahedral and dihedral solutions (of \cite{Dub95long, Hit-Poncelet, Hit-Octa}) fit into such families, as discussed in \cite{icosa, BH-Gav-p6families, cantat-loray07}. In general in such a family $y(t)$ may depend on the parameters of the family. Ben Hamed and Gavrilov \cite{BH-Gav-p6families} showed that any family with $y(t)$  {\em not} depending on the parameters is equivalent to one of the above cases and recently Cantat and Loray \cite{cantat-loray07} showed that any solution with $2,3$ or $4$  branches is in such family.

Next there is {\bf one discrete family} ($d,g$ unbounded, $\theta=(0,0,0,1)\sim (1,1,1,1)/2$). Indeed this \PVI\ equation was solved completely by Picard 
\cite{Picard89} p.299, R. Fuchs \cite{RFuchs} and in a different way by 
Hitchin \cite{Hit-tei}. Algebraic (determinantal)  formulae for the algebraic solutions amongst these appear in \cite{Hit-Poncelet}, using links with the Poncelet problem---in this framework they are dihedral solutions (controlling connections of type {\bf(A)} with binary dihedral monodromy).

Finally there are {\bf 45} {\bf exceptional solutions}, which collapse down to {\bf 30} if we identify solutions related by quadratic transformations. The possible genera are $0,1,2,3,7$ and the highest degree is $72$. Of these $30$ solutions $7$ have previously appeared: one is due to Dubrovin \cite{Dub95long}, two to Dubrovin--Mazzocco \cite{DubMaz00} and four to Kitaev (three in \cite{Kitaev-dessins}, plus---in \cite{kitaev-rmks}---a Galois conjugate of the elliptic $237$ solution already mentioned).  Two of these exceptional solutions are octahedral, one is the Klein solution, three are the elliptic $237$ solution (and its two Galois conjugates) and the remaining twenty-four are icosahedral.

\renewcommand{\baselinestretch}{0.9}              %
\normalsize
\bibliographystyle{amsplain}    \label{biby}
\bibliography{../thesis/syr} 

\vspace{0.5cm}   
\'Ecole Normale Sup\'erieure et CNRS, 
45 rue d'Ulm, 
75005 Paris, 
France

www.dma.ens.fr/$\sim$boalch

boalch@dma.ens.fr
\end{document}